\documentstyle[aps,pra,eqsecnum,12pt,preprint]{revtex}
\begin{document}
\author{D. Galetti, J.T. Lunardi and B.M. Pimentel}
\address{Instituto de F\'\i sica Te\'orica\\
Universidade Estadual Paulista\\
Rua Pamplona 145\\
01405 - 900 - S\~{a}o Paulo - SP\\
Brazil}
\author{C.L. Lima}
\address{Grupo de F\'{\i}sica Nuclear Te\'{o}rica e Fenomenologia de
Part\'{\i}culas\\
Elementares \\
Instituto de F\'{\i}sica }
\address{Universidade de S\~{a}o Paulo}
\address{Caixa Postal 66318}
\address{05389-970 S\~{a}o Paulo -- SP }
\address{Brazil}
\title{Unitary operator bases and q-deformed algebras}
\maketitle

\begin{abstract}
Starting from the Schwinger unitary operator bases formalism constructed out
of a finite dimensional state space, the well-known q-deformed commutation
relation is shown to emerge in a natural way, when the deformation parameter
is a root of unity.
\end{abstract}

\section{Introduction}

{}From the studies of deformed algebras, which appeared in connection with
problems in statistical mechanics and in quantum field theory (QFT), it came
out that the $q$-deformation parameter is, in its general form, a complex
number. In many applications it assumes a real value while in other cases
its imaginary part also plays a physical role. Apart from the basic quantum
mechanical study of the $q$-deformed oscillator by Biedenharn\cite{Bie} and
MacFarlane\cite{Mac}, in which a real deformation parameter is assumed,
Floratos\cite{Flo}, on the other hand, in his study of the $q$-oscillator
many-body problem, also discusses the case where $q$ is a pure complex
number. Furthermore,{\tt \ }in the particular case when $q$ is a root of
unity, it can be shown that the underlying state space, characterizing the
physical system, is finite dimensional. The $q$-deformed algebras generate a
suitable framework in this case and has been explicitly used in connection
with the phase problem in optics\cite{Abe}; moreover, it has also been
pointed out their importance in QFT\cite{Ber}.

A long time ago, Schwinger\cite{Sch} has pointed out that it is possible to
construct an operator basis, in the operator space, once we are given a
finite dimensional state space. The two fundamental unitary operators from
which the basis is constructed satisfy the Weyl commutation relation and act
cyclically on the corresponding state space, thus admitting as many roots of
unity, as eigenvalues, as is the dimension of the space. Here we will show
how the Schwinger operator basis can be used as a natural tool in order to
obtain the $q$-deformed commutation relation in the particular case when $q$
is a root of unity.

\section{The unitary operator bases}

For the complete quantum description of a physical system, a set of
operators must be found in such a way as to permit the construction of all
possible dynamical quantities related to that system. The elements of that
set are then identified as the elements of a complete operator basis.

One particular set, consisting of unitary operators, has been studied by
Schwinger \cite{Sch} and will be briefly recalled here. Let us consider a $N$
- dimensional linear, normed space of states to be understood as the quantum
phase-space of the relevant system. We can define a unitary operator $V$
through the mapping of an orthonormal system $\left\{ \langle u_k\mid
\right\} _{k=0...N-1}$ defined in this space, onto itself, by a cyclic
permutation as

\[
\langle u_k\mid V=\langle u_{k+1}\mid \text{ , }k=0,...N-1\text{, }\langle
u_N\mid =\langle u_0\mid \text{ .}
\]
A set of linearly independent unitary operators can then be constructed
trivially by mere repeated action of $V$,

\[
\langle u_k\mid V^s=\langle u_{k+s}\mid
\]
with

\[
\langle u_k\mid V^N=\langle u_k\mid \text{ ,}
\]
thus implying

\begin{equation}
V^N=\hat{1}\text{ , }  \label{2.1}
\end{equation}
where $\hat{1}$ is the unit operator.

The eigenvalues of $V$ obey this same equation and are thus given by the $N$
roots of unity

\[
v_k=\omega ^k=\exp \left( \frac{2\pi ik}N\right) .
\]
Furthermore, since that unitary operator has $N$ distinct eigenvalues, the
corresponding normalized eigenvectors, $\left\{ \langle v_l\mid \right\}
_{l=0,...N-1}$ , provide us with an alternative orthonormal system.

Schwinger has also shown that an operator $U$ exists such that

\[
\langle v_k\mid U=\langle v_{k-1}\mid \text{ ,}
\]
which is of period $N$, i.e.,

\[
U^N=\hat{1}\text{ , }
\]
thus implying the same spectrum as $V$ for the eigenvalues:

\[
u_k=\omega ^k=\exp \left( \frac{2\pi ik}N\right) .
\]
The fundamental point here is that the eigenvectors of $U$, $\left\{ \langle
u_k\mid \right\} _{k=0,...N-1}$ , can be shown to coincide with the
orthonormal set from which the construction started.

{}From (2.1) we can see that the special operator normalized to unit trace

\begin{equation}
\hat{G}(v_k)=\frac 1N\sum_{j=0}^{N-1}V^jv_k^{-j}  \label{2.2}
\end{equation}
is such that

\[
\langle v_l\mid \hat{G}(v_k)=\langle v_l\mid \delta _{l,k}\text{ ,}
\]
where

\[
\delta _{l,k}=\frac 1N\sum_{j=0}^{N-1}v_k^{-j}v_l^j
\]
plays the role of a Kronecker delta modulo $N$.

Corresponding to (2.2) we can also define

\begin{equation}
\hat{T}(u_k)=\frac 1N\sum_{j=0}^{N-1}U^{-j}u_k^j  \label{2.3}
\end{equation}
with additional equations similar to the above ones.

Using these properties we can show that the two coordinate systems are
related by a finite Fourier transformation with coefficients

\[
\langle u_k\mid v_l\rangle =\frac 1{\sqrt{N}}\exp \left( \frac{2\pi ikl}N%
\right) .
\]

Now, a simple verification leads us to the relation

\begin{equation}
V^lU^k=\exp \left( \frac{2\pi ikl}N\right) U^kV^l\text{ , }  \label{2.4}
\end{equation}
which, together with $V^N=\hat{1}$ and $U^N=\hat{1}$, fulfill the conditions
which characterize a generalized Clifford algebra \cite{Mo,Ya,Mor,Rama}.
Here, however, we will concentrate on just one special feature exhibited by
such a set of operators, viz., that the set of $N^2$ operators,

\[
\hat{S}_1\left( m,n\right) =\frac{U^mV^n}{\sqrt{N}}\text{ , }m,n=0,1,...,N-1%
\text{,}
\]
constitutes a complete orthonormal operators basis, with which we can
construct all possible dynamical quantities pertaining to that system \cite
{Sch}. In this way, an operator decomposition in this basis is written as

\begin{equation}
\hat{O}=\sum_{m,n=0}^{N-1}O\left( m,n\right) \hat{S}_1\left( m,n\right)
\text{ ,}  \label{2.5}
\end{equation}
where

\[
O\left( m,n\right) =Tr\left[ \hat{S}_1^{\dagger }\left( m,n\right) \hat{O}%
\right] .
\]

A very interesting property manifested by the operator basis $\left\{ \hat{S}%
_1\right\} $ is the factorization property

\[
\hat{S}_1\left( m,n\right) =\prod_{l=1}^h\hat{S}_{1l}\left( m_l,n_l\right)
\text{ , }
\]
where the sub-bases

\[
\hat{S}_{1l}\left( m_l,n_l\right) =\frac{U_l^{m_l}V_l^{n_l}}{\sqrt{P_l}}%
\text{ , }m_l,n_l=0,1,...P_l-1\text{ , }
\]
obey the commutation relations

\[
V_{l_1}U_{l_2}=U_{l_2}V_{l_1}\text{ , }l_1\neq l_2\text{ , }
\]

\[
V_{l_1}U_{l_2}=\exp \left( \frac{2\pi i}{P_{l_1}}\right) U_{l_2}V_{l_1}\text{
, }l_1=l_{2\text{ }},
\]
where $h$ is the total number of primes factors in $N$ including
repetitions, with $P_l$ a prime factor of $N$. This decomposition shows that
the factorized basis is constructed from operator sub-bases, each of which
associated with a prime number of states, the pair of operators $U$ and $V$
of each sub-basis being classified by the value of the prime integer $%
P_l=2,3,5,...$. It is straightforward to verify that the pair $U$ and $V$
associated with the canonical coordinate-momentum pair $q-p$ is obtained in
the particular case $P_l=\infty $. Then, according to Schwinger, due to this
factorization property and mutual orthogonality, each of these sub-bases is
associated to a particular degree of freedom of the physical system.

In order to emphasize the complete symmetry between $U$ and $V$, we want
also to observe that we could have introduced the new form for the operator
basis elements

\[
\hat{S}_2\left( m,n\right) =\frac{U^mV^n}{\sqrt{N}}\exp \left( \frac{i\pi mn}%
N\right) =\frac{V^nU^m}{\sqrt{N}}\exp \left( \frac{-i\pi mn}N\right) \text{
, }
\]
which preserves all properties already discussed under the substitutions $%
U\rightarrow V$ and $V\rightarrow U^{-1}$, combined with $m\rightarrow n$
and $n\rightarrow -m$.

For different degrees of freedom we must conveniently choose the range of
variation of the state labels in order to correctly treat the system
kinematics; for instance, it is important to emphasize again the canonical
case, i.e., $P_l=\infty $, for, in such a case, the unitary operators are
immediately identified with the well-known shift operators

\[
V\rightarrow e^{iq\hat{P}}
\]

\[
U\rightarrow e^{ip\hat{Q}}\text{ }
\]
when one considers the symmetric interval $m,n=-\frac{N-1}2,...,+\frac{N-1}2$%
, and then takes the $N\rightarrow \infty $ limit by prime numbers\cite{Sch}%
. However, it is also possible to perform a construction of the unitary
operators $U$ and $V$ in such a way to obtain an explicit ''angle - action''
pair, characterizing an Abelian two-dimensional rotation; in this case it
can be shown that

\begin{equation}
V\rightarrow \exp \left( i\frac{2\pi }N\hat{J}\right)  \label{2.6}
\end{equation}

\begin{equation}
U\rightarrow \exp \left( i\hat{\Theta}\right) .  \label{2.7}
\end{equation}

Here, the interval of variation of the state labels are suitably defined to
be $m=-\frac{N-1}2,...,+\frac{N-1}2$ and $n=-\frac{N-1}2\pi ,...,+\frac{N-1}2%
\pi $ in such a form that, in the limit of $N\rightarrow \infty $ , one
recovers $m=\left\{ -\infty ,...,+\infty \right\} $, running by integers,
and $n=\left\{ -\pi ,\pi \right\} $\cite{GaPi}.

For the sake of completeness it is important to go back to the operator
decomposition procedure, Eq.(2.5), and discuss the importance of the
particular choice of the operator basis. In fact, in order to emphasize the
discrete phase space character of the description, it was shown that\cite
{GaPi}, the Fourier transform of the original Schwinger basis $\hat{S}%
_2\left( m,n\right) $ must be considered so that a discrete Weyl transform
can be directly identified out of the decomposition scheme. With this new
basis it is straightforward to recover the well-known Weyl-Wigner
transformation for the canonical continuous case as well as a transformation
for an angle - angular momentum degree of freedom as special cases of the $%
N\rightarrow \infty $ limiting procedure. Furthermore, since the Schwinger
basis $\hat{S}_2\left( m,n\right) $ is not invariant under a modulo $N$
operation when the state labels are unrestricted in their domain, it was
shown that a new operator basis could be devised in order to preserve this
symmetry, namely\cite{GaPi2}

\[
\hat{G}\left( m,n\right) =\sum_{j=0}^{N-1}\sum_{l=0}^{N-1}\frac{\hat{T}%
\left( j,l\right) }{\sqrt{N}}\exp \left[ -\frac{2\pi i}N\left( mj+nl\right)
\right] \text{ , }
\]
where

\[
\hat{T}\left( j,l\right) =\hat{S}_2\left( j,l\right) \exp \left[ i\pi \phi
\left( j,l;N\right) \right] \text{ . }
\]
The phase $\phi \left( j,l;N\right) $ guarantees the mod $N$ invariance.

\section{ q-deformed algebras}

Since the Schwinger unitary operator bases formalism is constructed out of a
finite-dimensional state space, the relabelling procedure defines the
unitary shifting operators, which have as many eigenvalues (roots of unit)
as is the dimension of the underlying state space, $N$.

Let us now consider the set of eigenstates of the unitary operator $V$.
(Based on the symmetry stated in the last section, this choice is not
essential for what follows and could be replaced by the set of eigenstates
of the unitary operator $U$ as well.) Since $\left\{ \mid v_k\rangle
\right\} _{k=0,...N-1}$ is finite-dimensional and the unitary operator $U$
shifts cyclically the states of this space, one cannot interpret $U$ and $V$
as the corresponding creation and annihilation operators{\tt .} In fact, in
the space of eigenstates of $V$ the original pair of unitary operators are
represented as

\begin{eqnarray}
V &=&\sum_{l=0}^{N-1}\exp \left( \frac{2\pi il}N\right) \mid v_l\rangle
\langle v_l\mid =\sum_{l=0}^{N-1}v_l\mid v_l\rangle \langle v_l\mid
\nonumber  \label{3.1} \\
&=&\sum_{l=0}^{N-1}\omega ^l\mid v_l\rangle \langle v_l\mid  \label{3.1}
\end{eqnarray}

\begin{equation}
U=\sum_{l=0}^{N-1}\mid v_{l+1}\rangle \langle v_l\mid \text{ .}  \label{3.2}
\end{equation}

Nevertheless, starting from the unitary operators and making a convenient
choice for the state label range, we can construct a pair of operators which
will play the role of creation and annihilation operators in this
finite-dimensional space.

To begin with, it is immediate to see that, due to the symmetry of the
circle embodied in the unitary operator definition, one is not able to fix a
vacuum state solely from kinematical considerations, i.e., their action does
not select ''a priori'' any particular state as a vacuum state, since the
unitary operators act cyclically in the state space. In this case, one must
adopt some criterion to characterize this particular state. This choice will
break the symmetry of the circle and is not related to the kinematical
content of the description of the physical system. More precisely, one must
construct an operator out of the unitary operators in such a form to
annihilate the vacuum state; in addition, we must also have a creation
operator which generates the ''excited '' states of the multiplet.

The general form of the creation and annihilation operators will reveal the
possibility of particular choices for underlying algebras. To accomplish the
construction, we draw our attention again to the relations (3.1) and (3.2).
By comparison we can write these operators as

\[
a=\sum_{k=0}^{N-1}g\left( k\right) \mid v_k\rangle \langle v_{k+1}\mid
\]
and

\[
a^{\dagger }=\sum_{k=0}^{N-1}\mid v_{k+1}\rangle \langle v_k\mid .\text{ }
\]

The form of the creation operator only states that one jumps from the vacuum
state up to the last ($N-1$) state and so on cyclically. In what refers to
the annihilation operator we see that a particularly suitable choice for the
unknown function $g(k)$ must use an antisymmetric function of the state
label $k$ so as to select the vacuum state. Now, the natural antisymmetric
periodic function defined on the circle is the $\sin (\theta )$ function,
what requires odd $N$'s. Therefore, the proposed annihilation operator is
then written as

\[
a=\sum_{k=0}^{N-1}\frac{\sin \left( \frac{2\pi k}N\right) }{\sin \left(
\frac{2\pi }N\right) }\mid v_{k-1}\rangle \langle v_k\mid \text{ .}
\]

According to the discussion in section 2, we can decompose these operators
in the operator basis,

\[
\hat{O}=\sum_{m,n=0}^{N-1}O\left( m,n\right) \hat{G}\left( m,n\right)
\]
obtaining
\[
a^{\dagger }=U
\]
and

\[
a=U^{-1}\frac{V-V^{-1}}{\omega -\omega ^{-1}}
\]
respectively.

The question that can be posed now is if there exists some relation between
the bilinear products of the creation and annihilation operators, $%
a^{\dagger }$, $a$. Starting from the Weyl relation, Eq. (2.4), the
definitions Eqs. (2.6) and (2.7), where instead of $\hat{J}$ we now use $%
\hat{N}$ , the number operator and using the above expressions for $a$ and $%
a^{\dagger }$, we can immediately obtain the following relation

\[
aa^{\dagger }-\omega a^{\dagger }a=\omega ^{-\hat{N}}\text{ .}
\]

Therefore, we have seen that, starting from the Schwinger unitary operators,
the well-known $q$-deformed commutation relation emerges in a natural way,
when the deformation parameter is a root of unity. Furthermore, it is
immediate to verify that the creation and annihilation operators, $a$ and $%
a^{\dagger }$ proposed here are directly related to the ${\em h}$ and ${\em g%
}$ functions proposed by Floratos in his discussion of $q$-deformed algebras
for the bosonic case\cite{Flo}.

\section{Remarks and conclusions}

The main objective of this paper was to show that the $q$-deformed algebras
can be put in correspondence to Schwinger's unitary operator bases
formalism, when the deformation parameter is a root of unity.

Furthermore, it was shown that this formalism is the natural arena for the
discussion of recent work on general finite dimensional quantum mechanics
problems. Particularly, the Schwinger's formalism was used to represent any
operator acting on any finite dimensional state spaces \cite{GaPi,GaPi2}. To
be specific, it has also been used to study the Liouvillian dynamics in the
general finite dimensional phase spaces \cite{GaPi3} as well as to describe
physical systems from the particular case of a spin 1/2 ($N=2$ space) up to
the canonical continuous case (as the limit $N\rightarrow \infty $). The
special case of phase and number operators appearing in connection with
quantum optics has also been treated within this framework\cite{GaMa}. This
latter problem, or equivalently its Pegg-Barnett description\cite{Abe},
being just a particular case of the general Schwinger formalism, can
therefore be also embodied in the $q$-deformed algebra context along the
lines studied here.

\medskip\

{\bf Acknowledgment }D.G., J.T.L. and B.M.P. were supported by Conselho
Nacional de Desenvolvimento Cient\'{\i}fico e Tecnol\'{o}gico, CNPq, Brazil.

\end{document}